\documentclass[12pt]{amsart}
\usepackage{amssymb, amstext, amscd, amsmath, color, xypic}
\makeatletter
\def\@cite#1#2{{\m@th\upshape\bfseries%
[{#1\if@tempswa{\m@th\upshape\mdseries, #2}\fi}]}} \makeatother
\theoremstyle{plain}
\newtheorem{thm}[subsection]{Theorem}
\newtheorem{cor}[subsection]{Corollary}
\newtheorem{prop}[subsection]{Proposition}
\newtheorem{lem}[subsection]{Lemma}
\theoremstyle{definition}
\newtheorem{rem}[subsection]{Remark}

\newtheorem{defn}[subsection]{Definition}

\newtheorem{eg}[subsection]{Example}

\newcommand{\bC}{{\mathbb{C}}}

\newcommand{\bN}{{\mathbb{N}}}

\newcommand{\bT}{{\mathbb{T}}}

\newcommand{\bZ}{{\mathbb{Z}}}

\newcommand{\bm}{{\mathbf{m}}}
\newcommand{\bn}{{\mathbf{n}}}

  \newcommand{\A}{{\mathcal{A}}}

  \newcommand{\G}{{\mathcal{G}}}
\renewcommand{\H}{{\mathcal{H}}}
  
  \newcommand{\J}{{\mathcal{J}}}
  \newcommand{\K}{{\mathcal{K}}}

\renewcommand{\S}{{\mathcal{S}}}
  
  \newcommand{\U}{{\mathcal{U}}}
  
  \newcommand{\W}{{\mathcal{W}}}

\newcommand{\fA}{{\mathfrak{A}}}

\newcommand{\fF}{{\mathfrak{F}}}

\newcommand{\fM}{{\mathfrak{M}}}

\newcommand{\rC}{{\mathrm{C}}}

\newcommand{\ep}{\varepsilon}
\renewcommand{\phi}{\varphi}
\newcommand{\upchi}{{\raise.35ex\hbox{\ensuremath{\chi}}}}

\newcommand{\qand}{\quad\text{and}\quad}
\newcommand{\qif}{\quad\text{if}\quad}

\newcommand{\qfor}{\quad\text{for}\quad}
\newcommand{\qforal}{\quad\text{for all}\quad}

\newcommand{\AND}{\text{ and }}
\newcommand{\FOR}{\text{ for }}

\newcommand{\OR}{\text{ or }}

\newcommand{\id}{{\operatorname{id}}}

\newcommand{\ran}{\operatorname{Ran}}
\newcommand{\spn}{\operatorname{span}}

\newcommand{\Ath}{\A_\theta}

\newcommand{\ca}{\mathrm{C}^*}

\newcommand{\dlim}{\displaystyle\lim\limits}

\newcommand{\Fn}{\mathbb{F}_n^+}
\newcommand{\Fm}{\mathbb{F}_m^+}

\newcommand{\Fth}{\mathbb{F}_\theta^+}
\newcommand{\Ftheta}{\mathbb{F}_\theta^+}

\newcommand{\Fockth}{\ell^2(\Fth)}

\newcommand{\ip}[1]{\langle #1 \rangle}
\newcommand{\mt}{\varnothing}

\newcommand{\ol}{\overline}

\newcommand{\ltwo}{\ell^2}

\begin{document}
\title{Periodicity in Rank 2 Graph Algebras}

\author[K.R.Davidson]{Kenneth R. Davidson}
\address{Pure Math.\ Dept.\\U. Waterloo\\
Waterloo, ON\; N2L--3G1\\CANADA}
\email{krdavids@uwaterloo.ca}

\author[D.Yang]{Dilian Yang}
\address{Pure Math.\ Dept.\\U. Waterloo\\
Waterloo, ON\; N2L--3G1\\CANADA}
\email{dyang@uwaterloo.ca}

\begin{abstract}
Kumjian and Pask introduced an aperiodicity condition
for higher rank graphs.
We present a detailed analysis of when this occurs
in certain rank 2 graphs.
When the algebra is aperiodic, we give another proof
of the simplicity of $\ca(\Fth)$.
The periodic C*-algebras are characterized, and it is shown
that $\ca(\Fth) \simeq \rC(\bT) \otimes \fA$
where $\fA$ is a simple C*-algebra.
\end{abstract}

\thanks{First author partially supported by an NSERC grant.}

\date{}
\maketitle

\section{Introduction}\label{S:intro}
In this paper, we continue our study of the representation theory
of rank 2 graph algebras developed in \cite{DPYdiln} and \cite{DPYatomic}.
Kumjian and Pask \cite{KumPask} have introduced a family of
C*-algebras associated with higher rank graphs.
They describe a property called the aperiodicity condition
which implies the simplicity of the C*-algebra.
Our 2-graphs have a single vertex, and are particularly
amenable to analysis while exhibiting a wealth of
interesting phenomena.
Here we characterize when a 2-graph on one vertex is periodic,
and describe the associated C*-algebra.

The C*-algebras of higher rank graphs have been studied in a variety of papers
\cite{RaeSimYee1,RaeSimYee2,FarMuhYee,PaskRRS,RobSims,Sims}.
See also \cite{Raeburn}.
The corresponding nonself-adjoint algebras were introduced
by Kribs and Power \cite{KP1}.
The particular 2-graphs with one vertex were analyzed by Power \cite{P1}
and the representation theory was developed by Power and
us in \cite{DPYdiln, DPYatomic}.
Our work in this paper makes use of both sides of the theory; but this paper
is really about the C*-algebras.

In the case under consideration, the 2-graph is a semigroup $\Fth$
given by generators and relations.
We interpret the aperiodicity condition in terms of the existence of
a special faithful irreducible representation of the associated C*-algebra.
The typical situation is aperiodicity.
Indeed we show that periodicity only occurs under very special circumstances
in which the commutation relations for  words of certain lengths
are given by a flip operation.
Unfortunately, examples show that this periodicity may not exhibit
itself except for rather long words, making a determination in
specific examples difficult.
We develop an algorithm for doing the computations in a more
manageable way.

In the periodic case, there is also a special faithful representation.
It is not irreducible, but rather decomposes as a direct integral
by the methods of \cite{DPYatomic}.
The special structure in the periodic case allows us to provide a
detailed analysis of this direct integral, and thereby exhibit the
C*-algebra as a tensor product of $\rC(\bT)$ with a simple C*-algebra.
An important tool is a faithful approximately inner expectation onto
the C*-algebra generated by the gauge invariant AF-subalgebra and
the centre.

We would like to thank Stephen Power for helpful conversations.
We also thank Colin Davidson and  Shaoquan Jiang for assistance
in writing computer code for some of our calculations.

\section{Background}

The 2-graphs on a single vertex are semigroups which are given concretely
in terms of a finite set of generators and relations of a special type.
Let $\theta \in S_{m \times n}$ be a permutation of $\bm \times \bn$,
where $\bm=\{1,\dots,m\}$ and $\bn=\{1,\dots,n\}$.
The semigroup $\Fth$ is generated by $e_1,\dots,e_m$ and
$f_1,\dots,f_n$.
The identity is denoted as $\mt$.
There are no relations among the $e$'s, so they generate
a copy of the free semigroup on $m$ letters, $\Fm$;
and there are no relations on the $f$'s,
so they generate a copy of $\Fn$.
There are \textit{commutation relations} between the $e$'s
and $f$'s given by a permutation $\theta$ of $m \times n$:
\[ e_i f_j = f_{j'} e_{i'} \quad\text{where } \theta(i,j) = (i',j') .\]

A word $w\in\Fth$ has a fixed number of $e$'s and $f$'s regardless
of the factorization; and the \textit{degree} of $w$ is $(k,l)$ if
there are $k$ $e$'s and $l$ $f$'s.
The \textit{length} of $w$ is $|w| = k+l$.
The commutation relations allow any word $w\in\Fth$
to be written with all $e$'s first, or with all $f$'s first, say
$w = e_uf_v = f_{v'}e_{u'}$.
Indeed, one can factor $w$ with any prescribed pattern of
$e$'s and $f$'s as long as the degree is $(k,l)$.
It is straightforward to see that the factorization is
uniquely determined by the pattern and that $\Ftheta$
has the unique factorization property.
See also \cite{KumPask,P1}.

A representation $\sigma$ of $\Fth$ as operators on a Hilbert space
is \textit{row contractive} if
\[
 \big[ \sigma(e_1) \cdots \sigma(e_m)\big] \AND
 \big[\sigma(f_1) \cdots \sigma(f_n) \big]
\]
are row contractions; and \textit{row isometric}
if these row operators are isometries.
A row contractive representation is \textit{defect free} if
\[
 \sum_{i=1}^m \sigma(e_i) \sigma(e_i)^* = I
 = \sum_{j=1}^n \sigma(f_j) \sigma(f_j)^* .
\]
A row isometric defect free representation is called a
\textit{$*$-representation} of $\Fth$.
The universal C*-algebra for the family of $*$-representations
is denoted by $\ca(\Fth)$.
A faithful representation of $\ca(\Fth)$ will be denoted as $\pi_u$.

The left regular representation $\lambda$ of $\Fth$
is defined on $\Fockth$ with orthonormal basis
$\{ \xi_x : x \in \Fth\}$ by  $\lambda(w) \xi_x = \xi_{wx}$.
This is row isometric but is not defect free.
The norm closed unital operator algebra generated by these
operators is denoted by $\Ath$.

\subsection{Gauge automorphisms}
The universal property of $\ca(\Fth)$ yields a family of
\textit{gauge automorphisms}
$\gamma_{\alpha,\beta}$ for $\alpha, \beta \in \bT$
determined by
\[
 \gamma_{\alpha,\beta}(\pi_u(e_i)) = \alpha \pi_u(e_i) \qand
 \gamma_{\alpha,\beta}(\pi_u(f_j)) = \beta \pi_u(f_j) .
\]
Integration around the 2-torus yields a faithful expectation
\[ \Phi(X) = \int_{\bT^2}  \gamma_{\alpha,\beta}(X) \,d\alpha\,d\beta .\]
It is easy to check on monomials that the range is spanned
by words of degree $(0,0)$ (where $e_i^*$ and $f_j^*$ count as
degree $(-1,0)$ and $(0,-1)$ respectively).

Kumjian and Pask identify this range as an AF C*-algebra.
The first observation is that any monomial in $e$'s, $f$'s and
their adjoints can be written with all of the adjoints on the right.
Clearly the isometric condition means that
\[ \pi_u(f_i^*f_j) = \delta_{ij} = \pi_u(e_i^*e_j) .\]
To handle $e_i^*f_j$, observe that if
$f_j e_k = e_{k'}f_{j_k}$ for $1 \le k \le m$, then
\begin{align*}
 \pi_u(e_i^*f_j) &= \pi_u \big( e_i^*f_j \sum_k e_k e_k^* \big)\\
 &= \sum_k \pi_u( e_i^* e_{k'} f_{j_k} e_k^* )
 = \sum_k \delta_{ik'} \pi_u(  f_{j_k} e_k^* ) .
\end{align*}
So every word in $\ca(\Fth)$ can be expressed as a sum
of words of the form $xy^*$ for $x,y \in \Fth$.

Next, they observe that for each integer $k \ge1$, the set of words $\S_k$
in $\Fth$ of degree $(k,k)$ determine a family of degree $(0,0)$ words
$\{ \pi_u(xy^*) : x,y \in \S_k \}$.
It is clear that
\[ \pi_u(x_1y_1^*) \pi_u(x_2y_2^*) = \delta_{y_1,x_2} \pi_u(x_1y_2^*) .\]
Thus these operators form a family of matrix units that generate a
unital copy $\fF_k$ of the full matrix algebra $\fM_{(mn)^k}$.
Moreover, these algebras are nested because the identity
\[ \pi_u(xy^*) = \pi_u(x) \sum_i \pi_u(e_ie_i^*) \sum_j \pi_u(f_jf_j^*) \pi_u(y^*) \]
allows one to write elements of $\fF_k$ in terms of the basis for $\fF_{k+1}$.

Therefore the range of the expectation $\Phi$ is the
$(mn)^\infty$-UHF algebra $\fF = \ol{\bigcup_{k\ge1} \fF_k}$.
This is a simple C*-algebra.

\subsection{Type 3a representations}\label{3a_repn}
An important family of $*$-rep\-resent\-ations were introduced
in \cite{DPYdiln}.
The name refers to the classification obtained in \cite{DPYatomic}.

Start with an arbitrary \textit{tail} of $\Fth$, an infinite word of the form
\[ \tau = e_{i_0}f_{j_0}e_{i_1}f_{j_1} \dots .\]
Any infinite word $\tau$ with infinitely many $e$'s and infinitely many $f$'s
may be put into this \textit{standard form}.
It may also be factored with any pattern of $e$'s and $f$'s provided
that there are  infinitely many of each.
These alternate factorizations will be used later.

Let $\G_s = \G := \Ftheta$, for $s\ge0$, viewed as a discrete set
on which the generators of $\Ftheta$ act as injective maps
by right multiplication, namely,
\[ \rho(w)g = gw \qforal g \in \G. \]
Consider $\rho_s = \rho(e_{i_s}f_{j_s})$ as a map from
$\G_s$ into $\G_{s+1}$.
Define $\G_\tau$ to be the injective limit set
\[
 \G_\tau = \lim_{\rightarrow} (\G_s, \rho_s ) ;
\]
and let $\iota_s$ denote the injections of $\G_s$ into $\G_\tau$.
Thus $\G_\tau$ may be viewed as the union of $\G_0, \G_1, \dots $
with respect to these inclusions.

The left regular action $\lambda$ of $\Fth$ on itself induces
corresponding maps on $\G_s$ by  $\lambda_s(w) g = wg$.
Observe that $\rho_s \lambda_s(w) = \lambda_{s+1}(w) \rho_s$ .
The injective limit of these actions is an action $\lambda_\tau$
of $\Fth$ on $\G_\tau$.
Let $\lambda_\tau$ also denote the corresponding representation of $\Fth$
on $\ltwo(\G_\tau)$.
Let  $\{ \xi_g : g \in \G_\tau\}$ denote the basis.
A moment's reflection shows that this provides a defect free,
isometric representation of $\Fth$; i.e.\ it is a $*$-representation.
\medskip

It will be convenient to associate a directed chromatic graph to
any atomic representation $\sigma$.  We describe it for $\lambda_\tau$.
The vertices are associated to the points in $\G_\tau$.
For each vertex $x$ and each $i \in \bm$, draw a directed blue
edge labelled $e_i$ from $x$ to $y$ if $\lambda_\tau(e_i)\xi_x = \xi_y$.
Likewise for each $j \in \bn$, draw a directed red edge labelled
$f_j$ from $x$ to $z$ if $\lambda_\tau(f_j)\xi_x = \xi_z$.
Observe that defect free means that each vertex has one red and one blue
edge leading into the vertex.  For representations as partial isometries,
row contractivity means that there is at most one edge of each colour
leading into any vertex.  To be isometric, there must be $m$ blue edges and
$n$ red edges leading out of each vertex.

\medskip

One of the main results of \cite{DPYdiln} is that $\ca(\Fth)$ is the
C*-envelope of $\Ath$, and that every type 3a representation of $\Fth$
yields a completely isometric representation of $\Ath$ and a faithful
$*$-representation of $\ca(\Fth)$.

Therefore the gauge automorphisms are defined on $\ca(\lambda_\tau(\Fth))$.
It is shown in \cite{DPYdiln} that $\gamma_{\alpha,\beta}$ is implemented
on $\ltwo(\G_\tau)$ by the unitary
\[
 U_{\alpha,\beta} \xi_{\iota_s(e_uf_v)} =
 \alpha^{|u|-s}\beta^{|v|-s} \xi_{\iota_s(e_uf_v)} .
\]

\subsection{Coinvariant subspaces}
The other main result of \cite{DPYdiln} is that every defect free
representation of $\Fth$ extends to a completely contractive
representation of $\Ath$, and therefore dilates to a $*$-dilation.
Moreover the minimal dilation is unique.
Therefore it is possible to describe a $*$-representation completely
by its compression to a coinvariant cyclic subspace, as
it is then the unique $*$-dilation.

We describe such a subspace for the type 3a representations.
Let $\H = \ol{ \lambda_\tau(\Fth)^* \xi_{\iota_0(\mt)}}$.
This is coinvariant by construction.
As it contains $\xi_{\iota_s(\mt)}$ for all $s \ge1$,
it is easily seen to be cyclic.
Let $\sigma_\tau$ be the compression of $\lambda_\tau$ to $\H$.

Since $\lambda_\tau$ is a $*$-representation,
for each $(s,t) \in (-\bN)^2$, there is
a unique word $e_uf_v$ of degree $(|s|,|t|)$ such that
$\xi_{\iota_0(\mt)}$ is in the range of $\lambda_\tau(e_uf_v)$.
Set $\xi_{s,t} =  \lambda_\tau(e_uf_v)^* \xi_{\iota_0(\mt)}$.
It is not hard to see that this forms an orthonormal basis for $\H$.

Thus for each $(s,t) \in (-\bN)^2$, there are unique integers
$i_{s,t} \in \bm$ and $j_{s,t} \in \bn$ so that
\begin{alignat*}{2}
 \sigma_\tau( e_{i_{s,t}}) \xi_{s-1,t} &=\xi_{s,t} &\qfor&s\le 0 \AND t \le 0\\
 \sigma_\tau( f_{j_{s,t}}) \xi_{s,t-1} &=\xi_{s,t}  &\qfor&s \le 0 \AND t \le 0\\
 \sigma_\tau(e_i) \xi_{s,t} &= 0 &\qif& i \ne i_{s+1,t} \OR s = 0\\
 \sigma_\tau(f_j) \xi_{s,t} &= 0 &\qif& j \ne j_{s,t+1} \OR t = 0.
\end{alignat*}
Note that we label the edges leading \textit{into} each vertex,
rather than leading out, because in the $*$-dilation the blue (or red) edge
leading into a vertex is unique, while there are many leading out.

Consider how the  tail $\tau = e_{i_0}f_{j_0}e_{i_1}f_{j_1} \dots$
determines these integers.
It defines the path down the diagonal; that is,
\[ i_{s,s} := i_{|s|} \qand j_{s-1,s} :=j_{|s|} \qfor s \le 0 .\]
This determines the whole representation uniquely.
Indeed, for any vertex $\xi_{s,t}$ with $s,t \le 0$, take $T \ge |s|,|t|$,
and select a path from $(-T,-T)$ to $(0,0)$ that passes
through $(s,t)$.
The word $\tau_T = e_{i_0}f_{j_0} \dots e_{i_{T-1}}f_{j_{T-1}}$
satisfies $\sigma_\tau( \tau_T) \xi_{-T,-T} = \xi_{0,0}$.
Factor it as $\tau_T = w_1w_2$ with $d(w_1) = (T-|s|,T-|t|)$ and
$d(w_2) = (|s|,|t|)$, so that
$\sigma_\tau(w_2) \xi_{-T,-T} = \xi_{s,t}$
and $\sigma_\tau(w_1) \xi_{s,t} = \xi_{0,0}$.
Then $w_1 = e_{i_{s,t}}w' = f_{j_{s,t}} w''$.

It is evident that each $\sigma_\tau(e_i)$
and $\sigma_\tau(f_j)$ is a partial isometry.
Moreover, each basis vector is in the range of a unique
$\sigma_\tau(e_i)$ and $\sigma_\tau(f_j)$.
So this is a defect free, partially isometric representation
with unique minimal $*$-dilation $\lambda_\tau$.

\subsection{Symmetry and Periodicity}
An important part of the analysis of these atomic representations is the
recognition of symmetry.

\begin{defn}
The tail $\tau $ determines the integer data
\[ \Sigma(\tau) = \{ (i_{s,t},j_{s,t}) : s,t \le 0 \} .\]
Two infinite words $\tau_1$ and $\tau_2$ with data
$\Sigma(\tau_k) = \{ (i^{(k)}_{s,t},j^{(k)}_{s,t}) : s,t \le 0 \}$
are said to be \textit{tail equivalent} if the two sets of
integer data eventually coincide;
i.e.\ there is an integer $T$ so that
\[
 (i^{(1)}_{s,t},j^{(1)}_{s,t}) = (i^{(2)}_{s,t},j^{(2)}_{s,t})
 \qforal s,t \le T .
\]
Say that  $\tau_1$ and $\tau_2$ are \textit{$(p,q)$-shift tail equivalent}
for some $(p,q) \in \bZ^2$ if there is an integer $T$ so that
\begin{align}
 (i^{(1)}_{s+p,t+q},j^{(1)}_{s+p,t+q}) =
 (i^{(2)}_{s,t},j^{(2)}_{s,t}) \qforal s,t \le T . \tag{$*$}
\end{align}
Then $\tau_1$ and $\tau_2$ are \textit{shift tail equivalent} if they
are $(p,q)$-shift tail equivalent for some $(p,q) \in \bZ^2$.

The \textit{symmetry group} of $\tau$ is the subgroup of $\bZ^2$ given by
\[
 H_\tau = \{ (p,q) \in \bZ^2 :
 \Sigma(\tau) \text{ is $(p,q)$-shift tail equivalent to itself} \} .
\]
A sequence $\tau $ is called \textit{aperiodic} if $H_\tau = \{(0,0)\}$.

The semigroup $\Fth$ is said to satisfy the
\textit{aperiodicity condition} if there is an aperiodic infinite word.
Otherwise we say that $\Fth$ is \textit{periodic}.

We also say that $\tau$ is \textit{eventually $(p,q)$-periodic}
for $(p,q) \in H_\tau$.  If in fact it is fully $(p,q)$-periodic
(i.e. ($*$) holds whenever $s,t,s+p,t+q \le 0$), then we say that
$\tau$ is \textit{$(p,q)$-periodic}.
\end{defn}

In \cite{DPYatomic}, the atomic $*$-representations are
completely classified.
One of the important steps is defining a symmetry group
for the more general representations which occur.
It turns out that the representation is irreducible
precisely when the symmetry group is trivial.
So the aperiodicity condition is equivalent to saying that
there is an irreducible type 3a representation.

\section{Characterization of Periodicity}
\label{S:periodicity}

Whether or not there is an irreducible type 3a representation
of $\Fth$ depends on the semigroup.
In this section, we obtain detailed information about periodic
2-graphs.
In particular, a non-trivial symmetry group can only have the
form $\bZ(a,-b)$ where $a,b$ are integers such that $m^a=n^b$
and the commutation relations are very special.

Let $\bm^a$ denote the set of all $a$-tuples from the
alphabet $\bm$; and likewise $\bn^b$ denotes $b$-tuples
in the alphabet $\bn$.
We may (and shall) suppose that $m \le n$.
The case of $1=m < n$ is of limited interest, and $m=n=1$
is not considered.

\begin{thm}\label{T:periodic}
If $2 \le m \le n$, then the following are equivalent for $\Fth$
and positive integers $a$ and $b$.
\begin{enumerate}
\item Every tail of $\Fth$ is eventually $(a,-b)$ periodic.
\item Every tail of $\Fth$ is $(a,-b)$ periodic.
\item $m^a = n^b$ and there is a bijection $\gamma:\bm^a \to \bn^b$
so that
\[ e_u f_v = f_{\gamma(u)} e_{\gamma^{-1}(v)} \qforal u \in \bm^a \AND v \in \bn^b .\]
\end{enumerate}
If $1=m<n$, then $\Fth$ is $(0,b)$-periodic where $b$ is the order
of the permutation $\theta$.
\end{thm}

\begin{proof}
First consider $m \ge 2$, and
suppose that condition (3) holds.  Then we also have
$f_v e_u = e_{\gamma^{-1}(v)} f_{\gamma(u)}$.
Now consider a type 3a representation $\lambda_\tau$.
Fix any standard basis vector $\xi_{s,t}$ such that $s \le -a$.
Pulling back from $\xi_{s,t}$ yields a tail $\tau_1$,
which we factor as $\tau_1 = f_{v_1}e_{u_1}f_{v_2}e_{u_2} \dots$
where $|u_k|=a$ and $|v_k|=b$.
We wish to compare this with the tail $\tau_2$ obtained from $\xi_{s+a,t-b}$.

Note that starting at the vertex $\xi_{s+a,t}$, one gets to $\xi_{s,t}$ by
pulling back along a blue path $a$ steps using a word $e_u$; while one
obtains $\xi_{s+a,t-b}$ by pulling back $b$ steps along the red path $f_v$.
Hence the infinite path beginning at $\xi_{s+a,t}$ is
$\tau_0 = e_u \tau_1 = f_v \tau_2$.
Therefore
\begin{align*}
 \tau_0 &= e_u f_{v_1}e_{u_1}f_{v_2}e_{u_2} f_{v_3} \dots \\
    &= f_{\gamma(u)} e_{\gamma^{-1}(v_1)}  f_{\gamma(u_1)} e_{\gamma^{-1}(v_2)}
         f_{\gamma(u_2)} e_{\gamma^{-1}(v_3)} \dots \\
    &= f_v \tau_2 .
\end{align*}
Hence $v=\gamma(u)$ and
\begin{align*}
 \tau_2 &= e_{\gamma^{-1}(v_1)}  f_{\gamma(u_1)} e_{\gamma^{-1}(v_2)}
             f_{\gamma(u_2)} e_{\gamma^{-1}(v_3)} \dots \\
        &= f_{v_1}e_{u_1}f_{v_2}e_{u_2} f_{v_3} \dots = \tau_1.
\end{align*}
Therefore $\tau$ is $(a,-b)$-periodic.

Conversely, suppose that condition (3) fails.
We shall show that (1) also is false.
Condition (3) may fail for three reasons
relating to the identities $e_uf_v= f_{v'} e_{u'}$:
\begin{enumerate}
\renewcommand{\labelenumi}{(\roman{enumi}) }
\item $u'$ is not a function of $v$ alone,
\item $v'$ is not a function of $u$ alone, or
\item there are functions $\alpha: \bm^a \to \bn^b$ and $\beta:\bn^b \to \bm^a$
so that $e_uf_v= f_{\alpha(u)} e_{\beta(v)}$ but $\beta \ne \alpha^{-1}$.
\end{enumerate}

Consider (i) and select any $v\in \bn^b$ so that there are two words $u_i$
satisfying $e_{u_i} f_v= f_{v_i'} e_{u'_i}$ where $u'_1 \ne u'_2$.
Take an arbitrary word $u \in \bm^a$ and compute $f_v e_u = e_{u'} f_{v'}$.
Pick one of the $u_i$'s so that $u'_i \ne u'$.
Without loss of generality, this is $u_1$.
Now consider a word $e_{u_1}f_v e_u$ occurring as a segment of the
tail $\tau$, say $\tau = x e_{u_1}f_v e_u \tau'$.
In the 3a representation $\lambda_\tau$, there is a vertex $\xi_{s,t}$
at which the tail is
$\tau_0 = e_{u_1} f_v e_u \tau' = f_{v_1'} e_{u'_1} e_u \tau'$.
Moving to $\xi_{s-a,t}$ yields a vector with tail
$\tau_1 = f_v e_u \tau' = e_{u'} f_{v'} \tau'$.
Similarly, moving from $\xi_{s,t}$ to $\xi_{s, t-b}$ yields the tail
$\tau_2 = e_{u'_1} e_u \tau'$.
Since $u_1'\ne u'$, these two words do not coincide.

Hence any tail $\tau$ which contains the word $e_{u_1}f_v e_u$
infinitely often is not eventually $(a,-b)$ periodic.

Case (ii) is handled in the same manner.

In case (iii), note that this forces $\alpha$ and $\beta$ to be injections.
For $\alpha(u_1) = \alpha(u_2) = v_0$ implies that
$e_{u_1} f_v = f_{v_0} e_{\beta(v)} = e_{u_2} f_v$;
whence $e_{u_1} = e_{u_2}$ by cancellation.
Similarly for $\beta$.
Hence $m^a = n^b$ and $\alpha$ and $\beta$ are bijections.

Since $\beta \ne \alpha^{-1}$, select $v \in \bn^b$ so that
$\beta(v) \ne \alpha^{-1}(v)$.
Consider the tail $\tau = x e_{u_1} f_v e_{u_2} \tau'$.
Again there is a vertex $\xi_{s,t}$ at which the tail is
\[
 \tau_0 = e_{u_1} f_v e_{u_2} \tau' =
 f_{\alpha(u_1)} e_{\beta(v)} e_{u_2} \tau' .
\]
Moving to $\xi_{s-a,t}$ yields a vector with tail
$\tau_1 = f_v e_{u_2} \tau'= e_{\alpha^{-1}(v)} f_{\beta^{-1}(u_2)} \tau'$.
Similarly, moving from $\xi_{s,t}$ to $\xi_{s, t-b}$ yields the word
$\tau_2 = e_{\beta(v)} e_{u_2} \tau'$.
Since $\beta(v) \ne \alpha^{-1}(v)$, these two words do not coincide.
The proof is finished as before.

Now consider the case $m=1$.
Then $\theta \in S_n$ and the commutation relations have the form
$ef_j = f_{\theta(j)}e$ for $1 \le j \le n$.
So $e^k f_j = f_{\theta^k(j)} e$.
In particular, if $b$ is the order of $\theta$ in $S_n$,
then it is the smallest positive integer so that $e$
commutes with all $f_v$ for $v \in \bn^b$.

In this case, a type 3a representation is determined by the infinite
sequence $j_{0,t}$ for $t \le 0$.
Indeed, a simple calculation shows that $j_{s,t} = \theta^{-s}(j_{0,t})$
for all $s \le 0$.
Therefore every tail $\tau$ exhibits $(0,b)$ symmetry.
Select the sequence $(j_{0,t} : t \le 0)$ to be aperiodic
(as a sequence in one variable) and to contain all $n$ values
infinitely often.
It is easy to see that the data $\Sigma(\tau)$ exhibits
only $(0,b)$-periodicity.
\end{proof}

\begin{cor}
If $\frac{\log m}{\log n}$ is irrational, then
$\Fth$ is aperiodic for all $\theta$ in $S_{m \times n}$.
\end{cor}

\begin{proof}
$m^a=n^b$ if and only if $\frac{\log m}{\log n} = \frac b a$.
\end{proof}

\begin{eg}\label{periodic2x4}
Consider the following example with $m=2$ and $n=4$
with two 3-cycles (and two fixed points):
\[
 \big( (1,2), (2,1) , (1,3) \big) \qand
 \big( (2,2) , (2,3) , (1,4) \big).
\]
These relations are:
\begin{alignat*}{4}
e_1f_1&=f_1e_1 &\qquad e_1f_2 &= f_1e_2
&\qquad e_1f_3 &= f_2e_1 &\qquad e_1f_4 &= f_2e_2\\
e_2f_1&=f_3e_1 &\qquad e_2f_2 &= f_3e_2
&\qquad e_2f_3 &= f_4e_1 &\qquad e_2f_4 &= f_4e_2
\end{alignat*}
A calculation shows that the relation between $e$-words of length 2
and the $f$'s has this special symmetry.  Setting
\[ \gamma(11) = 1 \qquad \gamma(12) = 2 \qquad \gamma(21) = 3 \qquad \gamma(22) = 4 \]
yields the relations $e_{ij} f_k = f_{\gamma(ij)} e_{\gamma^{-1}(k)}$.
So this semigroup has $(2,-1)$-periodicity.
\end{eg}

Theorem~\ref{T:periodic} leads to the following theorem.
It is somewhat unsatisfactory because one needs to check possibly
many higher commutation relations (see Example~\ref{twelve}).
We pose the question of whether there is a
combinatorial condition on the original permutation $\theta$
which is equivalent to periodicity.
A partial answer to this problem is given later in this section.

\begin{thm}\label{aperiodic}
Suppose that $m,n \ge 2$.  Then
$\Fth$ satisfies the aperiodicity condition if and only if
the technical condition $(3)$ of Theorem~$\ref{T:periodic}$
does not hold for any $(a,b)$ for which $m^a=n^b$.
\end{thm}

\begin{proof}
List all non-zero words $(p,q) \in \bZ^2$ in a list $\{(p_k,q_k) : k \ge1\}$
so that each element is repeated infinitely often.
For each $k$, we construct a word $a_k$ in $\Fth$.
There are two cases for the word $(p,q)$.

If $p_kq_k\ge0$ and $p_k\ne0$, choose $a_k = e_1^{|p_k|} f_1^{|q_k|} e_2$;
and if $p_k=0$, choose $a_k=f_1^{|q_k|}f_2$.
If $p_kq_k<0$, use the construction from the proof of Theorem~\ref{T:periodic}.
Let $\tau = a_1a_2a_3\dots$.

To see that $\tau$ is aperiodic, consider any $(p,q) \ne (0,0)$.
It occurs as $(p_k,q_k)$ infinitely many times.
If $p_kq_k\ge0$, consider the starting point $(s,t)$ at the beginning
of the word $a_k= e_1^{|p_k|} f_1^{|q_k|} t$ where $t=e_2$ or $f_2$.
Then moving to $(s-|p_k|,t-|q_k|)$ yields the word beginning with $t$,
which does not coincide with the beginning of $a_k$.
If $p_kq_k<0$, then argue as in the previous Theorem.
As each $(p,q)$ occurs infinitely often, $\tau$ is not eventually
$(p,q)$-periodic for any period.
Hence it is aperiodic.
\end{proof}

The same proof works for the periodic semigroups, eliminating all
symmetries except those in every representation.

\begin{cor}\label{C:periodic}
If $\Fth$ is periodic with minimal period $(a,-b)$,
then $\ca(\Fth)$ has a type 3a representation with symmetry
group $\bZ(a,-b)$.
\end{cor}

\begin{proof}
If $\Fth$ is periodic with minimal period $(a,-b)$, then a routine modification
of the proof above shows that there is an infinite word whose only symmetries
are $\bZ(a,-b)$.
Indeed, if $m^{a_0}=n^{b_0}$ and $\gcd(a_0,b_0)=1$, then
$a=ka_0$ and $b=kb_0$.
By hypothesis, there are words $a_l$ with no $(la_0,-lb_0)$
periodicity for $1 \le l < k$.
As in the proof above, there are always words with no $(p,q)$ periodicity
when $pq\ge0$.
So following the same process, one obtains an infinite word $\tau$
without any of these symmetries.
Let $\lambda_\tau$ be the corresponding type 3a representation.

The symmetry group $H_{\lambda_\tau}$ contains $(a,-b)$,
and thus $\bZ(a,-b)$.
However by construction, $H_{\lambda_\tau} \cap \bN_0^2 = \{(0,0)\}$ and
$(la_0,-lb_0)$ are not in $H_{\lambda_\tau}$ for $1 \le l < k$.
So $H_{\lambda_\tau} \cap \bZ(a_0,-b_0) = \bZ(a,-b)$.
If it contains anything else, then it would contain a non-zero
element of $\bN_0^2$.
So $H_{\lambda_\tau} = \bZ(a,-b)$.
\end{proof}

\section{Tests for Periodicity and Aperiodicity}
\label{S:tests}

We now examine a method for demonstrating aperiodicity.
The permutation $\theta \in S_{mn}$ determines functions
$\alpha_i : \bn \to \bn$ and $\beta_j : \bm \to \bm$ so that
\[ \theta(i,j) = (\beta_j(i), \alpha_i(j)) .\]
Thus if $u=i_1\dots i_a$ and $v = j_b\dots j_1$,
\[ e_u f_j = f_{\alpha_{i_1} \circ \alpha_{i_2} \circ \dots \circ \alpha_{i_a}(j)} e_{u'} \]
and
\[ e_i f_v = f_{v'} e_{\beta_{j_1} \circ \beta_{j_2} \circ \dots \circ \beta_{j_b}(i)} .\]

For $(a,-b)$-periodicity when $m^a=n^b$, a necessary condition is that
$\alpha_{i_1} \circ \alpha_{i_2} \circ \dots \circ \alpha_{i_a}$ and
$\beta_{j_1} \circ \beta_{j_2} \circ \dots \circ \beta_{j_b}$
are constant maps for all $u$ and $v$.
So we obtain:

\begin{cor}[Aperiodicity criterion]
If there is a subset $B \subset \bn$ with $|B|\ge2$
and a word $i_1\dots i_k \in \bm^k$ so that
$\alpha_{i_1} \circ \alpha_{i_2} \circ \dots \circ \alpha_{i_k}(B)=B$,
then $\Fth$ is aperiodic.
Similarly, if there is a subset $A \subset \bm$ with $|A|\ge2$
and a word $j_1\dots j_k \in \bn^k$ so that
$\beta_{j_1} \circ \beta_{j_2} \circ \dots \circ \beta_{j_k}(A)=A$,
then $\Fth$ is aperiodic.
\end{cor}

\begin{proof}
$\big( \alpha_{i_1} \circ \alpha_{i_2} \circ \dots \circ
\alpha_{i_k} \big)^a(B) = B$ is never constant.
\end{proof}

\begin{rem}
It is not hard to show that either there is a
$B \subset \bn$ with $|B|\ge2$
and a word $i_1\dots i_k \in \bm^k$ so that
$\alpha_{i_1} \circ \alpha_{i_2} \circ \dots \circ \alpha_{i_k}(B)=B$
or for some sufficiently large $k$, all
$\alpha_{i_1} \circ \alpha_{i_2} \circ \dots \circ \alpha_{i_k}$
are constant maps.
\end{rem}

\begin{eg}\label{2by2}
For $\theta \in S_{2\times2}$, there are 9 distinct algebras up
to isomorphism \cite{P1}.

For example, the forward 3-cycle algebra is given by
the permutation $\theta$ in $S_{2 \times 2}$ given by
the $3$-cycle  $\big( (1,1) , (1,2) , (2,1) \big)$.
This yields the relations
\begin{alignat*}{2}
e_1f_1 &= f_2e_1, &\qquad e_1f_2 &= f_1e_2 \\
e_2f_1 &= f_1e_1, &\qquad e_2f_2 &= f_2e_2.
\end{alignat*}
One can easily check that $e_if_j = f_{i+j}e_j$,
where addition is modulo 2.
Notice that $\alpha_2 = \id$;
so $\alpha_2 (\{1,2\}) = \{1,2\}$.
Hence it is aperiodic.
This technique works for seven of the nine $2 \times 2$ examples.

One exception is the flip algebra, which is
given by the rule $e_if_j = f_ie_j$;
and it is clearly $(1,-1)$-periodic.

The other is the square algebra given by the permutation
\[ \big( (1,1) , (1,2) , (2,2) , (2,1) \big) .\]
This yields maps
\[ \alpha_1 = \beta_2 = 2  \qand  \alpha_2 = \beta_1 = 1 .\]
These maps are constant, but are not mutual inverses.
So there is no $(1,-1)$ periodicity.
However a calculation shows that
\[ e_{i_1i_2} f_{j_1j_2} = f_{i'_1i_2} e_{j'_1j_2} \]
where $i' = i+1 \pmod2$ and $j'=j+1\pmod2$.
So the function $\gamma(ij)=i'j$ satisfies $\gamma^{-1}=\gamma$
and
\[
  e_uf_v = f_{\gamma(u)} e_{\gamma(v)} \qforal |u|=|v|=2 .
\]
Thus the square algebra is $(2,-2)$-periodic.

So the periodicity can reveal itself only in the higher order
commutation relations!
\end{eg}

\begin{eg}\label{periodic3x3}
Here is another example of this phenomenon.
Consider $\theta\in S_{3\times 3}$ given by fixed points 
$(i,i)$ for $1 \le i \le 3$, and cycles
\[ \big( (1,2), (2,1)\big) \qand \big( (1,3) , (3,2) , (2,3) , (3,1) \big) . \]
A calculation shows that this algebra has $(2,-2)$-periodicity via
the correspondence $e_{ij} f_{kl} = f_{\gamma(ij)} e_{\gamma(kl)}$
where $\gamma(23) = 13$, $\gamma(13)=23$ and $\gamma(ij)=ij$ otherwise
(so $\gamma^{-1}=\gamma$).

These calculations can be simplified somewhat by observing that
there are subalgebras isomorphic to the one in Example~\ref{periodic2x4}
generated by each of the sets
$\{e_1,e_2; f_{11}, f_{12}, f_{21}, f_{22}\}$,
$\{e_1,e_3; f_{11}, f_{23}, f_{31}, f_{33}\}$,
$\{e_2,e_3; f_{22}, f_{13}, f_{32}, f_{33}\}$,\!
$\{f_1,f_2; e_{11}, e_{12}, e_{21}, e_{22}\}$,\!
$\{f_2, f_3; e_{22}, e_{13}, e_{32}, e_{33}\}$, and
$\{f_1, f_3; e_{11}, e_{23}, e_{31}, e_{33}\}$.
Thus there are corresponding $4\times 4$ subsets of the $9 \times 9$
pattern of relations between words of degree 2 that must have the
desired form.
\end{eg}

In order to develop a better test, we require a refinement of 
condition (3) of Theorem~\ref{T:periodic}.

\begin{prop}\label{periodic identity}
If $\Fth$ is periodic and $\gamma: \bm^a\to \bn^b$
is the bijective correspondence of Theorem~$\ref{T:periodic}$,
then for $i_0,\dots,i_a \in \bm$,
\[
 e_{i_0} f_{\gamma(i_1 \dots i_a)} =  f_{\gamma(i_0 \dots i_{a-1})} e_{i_a} .
\]
Conversely, if there is a bijection $\gamma: \bm^a \to \bn^b$ with
this property, then $\Fth$ is periodic.

Similarly for $j_0,\dots,j_b \in \bn$,
$
 e_{\gamma^{-1}(j_0 \dots j_{b-1})} f_{j_b} =  f_{j_0} e_{\gamma^{-1}(j_1 \dots j_b)}
.$ 
Again this property is equivalent to $(a,-b)$-periodicity.
\end{prop}

\begin{proof}
From the commutation relations, we know that
\[
 e_{i_0} f_{\gamma(i_1 \dots i_a)} = f_v e_{i'_a} =   f_{\gamma(i'_0 \dots i'_{a-1})} e_{i'_a}
\]
where  $i'_0 \dots i'_{a-1} = \gamma^{-1}(v)$.
Let $u = ku'$ be any word of length $a$.
Then
\begin{align*}
 e_u e_{i_0} f_{\gamma(i_1 \dots i_a)} &=
 e_k e_{u'i_0} f_{\gamma(i_1 \dots i_a)}  =
 e_k f_{\gamma(u'i_0)} e_{i_1 \dots i_a} \\ &=
 e_{u} f_{\gamma(i'_0 \dots i'_{a-1})} e_{i'_a} =
 f_{\gamma(u)} e_{i'_0 \dots i'_{a-1}i'_a} .
\end{align*}
Therefore $i'_s=i_s$ for $1 \le s \le a$.
Similarly,
\begin{align*}
 e_{i_0} f_{\gamma(i_1 \dots i_a)} e_{u'} &=
 f_{\gamma(i'_0 \dots i'_{a-1})} e_{i'_au'}   =
 e_{i'_0 \dots i'_{a-1}} f_{\gamma(i'_au')} .
\end{align*}
So $i'_0=i_0$.

Conversely, using this identity $a$ times yields
\[ e_u f_{\gamma(u')} = f_{\gamma(u)} e_{u'} \qforal u,u' \in \bm^a .\]
Thus $\Fth$ is periodic by Theorem~\ref{T:periodic}.
\end{proof}

The next proposition shows that the periodicity of the square algebra
of Example~\ref{2by2} is a consequence of the relation
$e_if_j = f_{i+1}e_j$.

\begin{prop}\label{alpha beta periodic}
If  $m^a=n^b$ and there are maps $\alpha:\bm^a \to \bn^b$
and $\beta:\bn^b \to \bm^a$ such that $e_uf_v = f_{\alpha(u)} e_{\beta(v)}$
for all $u \in \bm^a$ and $v \in \bn^b$, then $\Fth$ is periodic.
\end{prop}

\begin{proof}
Since there is a bijective correspondence $\theta'$ between
the words $e_uf_v$ and the words $f_{v'}e_{u'}$, it is easy to
verify that $\alpha$ and $\beta$ are bijections.
Thus $\beta\alpha$ is a permutation of $\bm^a$.
Let $k$ be the order of $\beta\alpha$ in $S_{m^a}$;
so $(\beta\alpha)^k = \id$.
Define $\gamma : \bm^{ak} \to \bn^{bk}$ by
\[
 \gamma(u_1 \dots u_k) =
 \alpha(u_1) \alpha\beta\alpha(u_2) \dots \alpha(\beta\alpha)^{k-1}(u_k)
\]
for $u_i \in \bm^a$, $1 \le i \le k$.
Then compute
\begin{align*}
 e_{u_0} f_{\gamma(u_1 \dots u_k)} &=
 e_{u_0} f_{\alpha(u_1)} f_{\alpha\beta\alpha(u_2)}
  \dots f_{\alpha(\beta\alpha)^{k-1}(u_k)} \\ &=
 f_{\alpha(u_0)} e_{\beta\alpha(u_1)} f_{\alpha\beta\alpha(u_2)}
  \dots f_{\alpha(\beta\alpha)^{k-1}(u_k)} \\ &=
 f_{\alpha(u_0)} f_{\alpha\beta\alpha(u_1)} \dots
 f_{\alpha(\beta\alpha)^{k-1}(u_{k-1})} e_{(\beta\alpha)^k(u_k)} \\&=
 f_{\gamma(u_0u_1\dots u_{k-1})} e_{u_k} .
\end{align*}
Therefore by Proposition~\ref{periodic identity},
$\Fth$ is $(ak,-bk)$-periodic.
\end{proof}

Proposition~\ref{periodic identity} can be used to calculate $\gamma$.
In turn, this leads to a checkable algorithm for periodicity.
This is captured in the following theorem.

\begin{thm}\label{T:algorithm}
Suppose that $\Fth$ is $(a,-b)$-periodic.
Then the bijection $\gamma: \bm^a \to \bn^b$ may be
calculated for $u_0 \in \bm^a$ by starting with an
arbitrary $j_0 \in \bn$ and computing
\begin{align*}
 e_{u_0} f_{j_0} &= f_{j_1} e_{u_1}\\
 e_{u_1} f_{j_1} &= f_{j_2} e_{u_2}\\
            &\vdots \\
 e_{u_b} f_{j_b} &= f_j e_u
\end{align*}
Then $\gamma(u_0) =  v_0 := j_1j_2\dots j_b$,
and also $j=j_0$ and $u=u_0$.

Reversing the process, start with an arbitrary
$i_0 \in \bm$ and calculate
\begin{align*}
 e_{i_0} f_{v_0} &= f_{v_1} e_{i_1}\\
 e_{i_1} f_{v_1} &=  f_{v_2} e_{i_2}\\
            &\vdots \\
 e_{i_a} f_{v_a} &= f_v e_i.
\end{align*}
Then $\gamma^{-1}(v_0) = i_ai_{a-1} \dots i_1 = u_0$,
and also $v=v_0$ and $i=i_0$.

Conversely, if for $m^a=n^b$ and for each $u_0 \in \bm^a$,
the procedure above passes all the tests of equality
for all $j_0 \in \bn$ and all $i_0\in\bm$, then
$\Fth$ is $(a,-b)$-periodic.
\end{thm}

\begin{proof}
By Proposition~\ref{periodic identity}, if $\gamma(u_0) = j_1j_2 \dots j_b$,
then
\[
 e_{u_0} f_{j_0} = f_{j_1} e_{u_1} \quad\text{where}\quad
 u_1 = j_2\dots j_b j_0 .
\]
Proceeding by induction, we find that
$u_i = j_{i+1} \dots j_b j_0 \dots j_{i-1}$ and
\[ e_{u_i} f_{j_i} = f_{j_{i+1}} e_{u_{i+1}} \]
where $u_{i+1} = j_{i+2} \dots j_b j_0 \dots j_i$.
In the last step, we return to the beginning and obtain
\[ e_{u_b} f_{j_b}  = f_{j_0} e_{u_0} .\]
Hence we have calculated $\gamma(u_0) = j_1j_2 \dots j_b$,
and $j=j_0$ and $u=u_0$.

Reversing the process works in the same manner.

Now consider the converse.
Starting with each $u_0$ in $\bm^a$, for each value of $j_0 \in \bn$,
we produce the same sequence $v_0 := j_1j_2\dots j_b$.
This defines a function $\alpha:\bm^a \to \bn^b$.
Observe that with the notation from that calculation,
since the sequence cycles around due to the fact that $v=v_0$ and $j=j_0$,
it follows that $\alpha(u_i) = j_{i+1} \dots j_b j_0 \dots j_{i-1}$
for $1 \le i \le b$.

Then we reverse the process, and construct a function
$\beta:\bn^b \to \bm^a$ and confirm that $\beta(v_0) = u_0$.
That is, $\beta = \alpha^{-1}$.
Therefore $\alpha$ and $\beta$ are bijections.
Finally, the initial calculation $e_{u_0} f_{j_0} = f_{j_1} e_{u_1}$
yields
\[ e_{\alpha^{-1}(j_1 \dots j_b)} f_{j_0} = f_{j_1} e_{\alpha^{-1}(j_2 \dots j_b j_0)} .\]
This verifies the hypothesis of Proposition~\ref{periodic identity},
and confirms that $\Fth$ is $(a,-b)$-periodic.
\end{proof}

\subsection{A computer algorithm}
This theorem provides a valid test for periodicity that is effective as a computer program.
It allows a single pass through all the words in $\bm^a$ doing several tests.
Failure at any point indicates failure of $(a,-b)$-periodicity; while a completed run without
failure means that $\Fth$ is indeed  $(a,-b)$-periodic.
An algorithm based on Theorem~\ref{T:periodic} would require checking all
$m^an^b = m^{2a}$ pairs, and this is much too computationally intensive.

\begin{eg}\label{twelve}
This example is a $4\times4$ example which has $(12,-12)$-periodicity.
This is surprisingly high periodicity for such a small number of
generators.  Already it is basically impossible to calculate the
multiplication table for the $4^{12} \times 4^{12}$ pairs of words.
The algorithm described above reduces this example to a calculation
that can be done by computer in about an hour.
We first show how hand calculations allow us to deduce that
there is no periodicity smaller than $12$.

We call this the 8-cycle algebra.  It is given by the 8-cycle:
\[ \big( (2,1), (1,2) , (3,1) , (1,3) , (4,2) , (2,4) , (4,3) , (3,4) \big), \]
two 2-cycles $\big( (1,4), (4,1) \big)$ and $\big( (2,3) , (3,2) \big)$,
and fixed points $(i,i)$ for $1 \le i \le 4$.

We first calculate the maps $\alpha_i$ and $\beta_j$.
\begin{alignat*}{2}
\alpha_1(3) &= 2 &\qand \alpha_1(j)&=1 \text{ otherwise} \\
\alpha_2(4) &= 3 &\qand \alpha_2(j)&=2 \text{ otherwise} \\
\alpha_3(4) &= 1 &\qand \alpha_3(j) &= 3 \text{ otherwise} \\
& &\quad \alpha_4 &= 4\\
& &\quad \beta_1 &=1  \\
\beta_2(1)&=3 &\qand \beta_2(i) &= 2 \text{ otherwise} \\
\beta_3(1) &= 4 &\qand \beta_3(i) &= 3 \text{ otherwise} \\
\beta_4(3) &= 2 &\qand \beta_4(i) &= 4 \text{ otherwise}
\end{alignat*}
One readily calculates
$\alpha_2^2 = \alpha_2\alpha_1 = \alpha_2\alpha_3 = 2$,
$\alpha_3^2 = \alpha_3\alpha_1 = \alpha_3 \alpha_2 = 3$,
 $\alpha_1^2=1$,
and any expression involving $\alpha_4$ is constant as well.
However $\alpha_1\alpha_3$ and $\alpha_1\alpha_2$
are not constant.  It follows though that every composition
of three $\alpha_i$'s is constant.
A similar calculation shows that the composition of
any three $\beta_j$'s is constant.
This suggests that the 8-cycle algebra may be periodic.
But one might think that it should have small order.
That turns out not to be the case.

We will show by hand that if the 8-cycle algebra $\Fth$ is
periodic with the minimal period $(a,-b)$, 
then $(a,-b)=(12k,-12k)$ for some $k\ge 1$.
Clearly $m=n$ implies that $a=b$.
A useful observation is that if $\Fth$ is $(k,-k)$-periodic
then when $|u|=|v|=k$ and $e_uf_v=f_{v'}e_{u'}$, 
it follows that $e_{u'}f_{v'}=f_ve_u$.
That is, the cycle lengths are just 1 and 2.
We will show that this forces $k$ to be a multiple of 12.

Observe that in the commutation relations 
between the 2 letter words $\{11,12,13,24,34\}$
remain within this set; and so we obtain a subsemigroup
generated by these words that is a 2-graph with
a $5\times 5$ multiplication table.

The point $(12, 11)$ lies on the 6 cycle
\[ \big( (12, 11),(11,13),(34,11),(11,12),(13,11),(11,24) \big) .\]
By induction, we can obtain the following identities
\begin{alignat*}{2}
e_{12}^{2k+1} f_{11}^{2k+1} &= (f_{13}f_{34})^k f_{13} e_{11}^{2k+1} 
&\qfor &k \ge 0\\
e_{12}^{2k} f_{11}^{2k} &= (f_{13}f_{34})^k e_{11}^{2k} 
&\qfor &k \ge 1 .
\end{alignat*}
On the other hand, we compute
{\allowdisplaybreaks
\begin{align*}
e_{11}^{6k+1} (f_{13}f_{34})^{3k} f_{13}   &= f_{11}^{6k+1} e_{34} e_{24}^{6k} \\
e_{11}^{6k+2} (f_{13}f_{34})^{3k+1}        &= f_{11}^{6k+2} e_{13} e_{12}^{6k+1} \\
e_{11}^{6k+3} (f_{13}f_{34})^{3k+1} f_{13} &= f_{11}^{6k+3} e_{12} e_{24}^{6k+2} \\
e_{11}^{6k+4} (f_{13}f_{34})^{3k+2}        &= f_{11}^{6k+4} e_{34} e_{12}^{6k+3} \\
e_{11}^{6k+5} (f_{13}f_{34})^{3k+2} f_{13} &= f_{11}^{6k+5} e_{13} e_{24}^{6k+4} \\
e_{11}^{6k+6} (f_{13}f_{34})^{3k+3}        &= f_{11}^{6k+6} e_{12}^{6k+6}.
\end{align*}
} 

 From this, we see that the required 2-cycle condition does not
hold for words of length $12k+2i$ for $1 \le i \le 5$.
It also follows that there is no odd period $(k,-k)$,
for then $(2k,-2k)$ would be a period which is not a multiple of 12.
Therefore, $(k,-k)$-periodicity can only hold if $12$ divides $k$.

The computer algorithm successfully verified
that $\Fth$ is $(12,-12)$-periodic.
Hence the symmetry group is $\bZ(12,-12)$.

As a corollary, we see that the $5\times5$ algebra that we used 
has symmetry group $\bZ(6,-6)$.
This follows because the map $\gamma$ on $4^{12}$ restricts
to a map on the 6-letter words from the $5\times5$ algebra.
So there is $(6,-6)$-periodicity.
Our argument shows that it has no smaller period.
\end{eg}

For a while, we had conjectured that if there is a constant $k$
so that the composition of any $k$ of the maps $\alpha_i$ is constant,
as is the composition of any $k$ of the maps $\beta_j$, then $\Fth$ 
should be periodic.
However the following example shows that this is not the case.
So we pose the less precise problem: \textit{
Find a computable condition on the permutation $\theta$
which is equivalent to periodicity.}

\begin{eg}
Consider the $3\times 3$ example with an $8$-cycle:
\[ \big( (1,3),(1,1),(3,1),(3,3),(2,3),(1,2),(2,1),(3,2) \big) .\]
It is easy to calculate that $\alpha_i = i$ are constant
for $i=1,2,3$, as is $\beta_1 = 3$; $\beta_2$ sends $3$ to $1$
and $1,2$ to $2$; and $\beta_3$ sends $3$ to $2$, and $1,2$ to $1$.
So $\beta_2^2 = 2 = \beta_2\beta_3$ and $\beta_3^2 = 1 = \beta_3\beta_2$.
So all compositions of two maps are constant.

We claim that $e_{1^k} f_{1^k} = f_{132^{k-2}} e_{u_k}$.
Indeed this is an easy calculation by induction starting with
$e_{11}f_{11}=f_{13}e_{23}$, since
\begin{align*}
 e_{1^{k+1}} f_{1^{k+1}} &= e_1 f_{132^{k-2}} e_{u_k} f_1 
 = f_{13} e_2 f_{2^{k-2}} f_j e_{u'_k} \\
 &= f_{132^{k-2}} e_2 f_j e_{u'_k} = f_{132^{k-2}} f_2 e_{i'u'_k} = f_{132^{k-1}} e_{u_{k+1}} .
\end{align*} 
Thus if $\Fth$ is $(k,-k)$-periodic for $k \ge4$,
one would have to have $e_{u_k} f_{132^{k-2}} = f_{1^k}e_{1^k}$.
However
\[
 e_{u_k} f_{132^{k-2}} = 
 e_{u'_k} e_i f_{132^{k-2}} = 
 e_{u'_k} f_{j_1j_2} e_{i'} f_{2^{k-2}} =
 e_{u'_k} f_{v'} e_2 
\]
because $\beta_2^2 = 2$.  Thus $\Fth$ must be aperiodic.
\end{eg}

\begin{eg}
Here is another example with a $4\times 4$ permutation:
\begin{alignat*}{2}
 &\big((1,1), (3,2), (4,4), (2,3) \big),
 &\qquad&\big((2,1), (1,2), (4,2), (2,4) \big),\\
 &\big((3,1), (3,4), (4,3), (1,3) \big)
 &\qquad&\big ((4,1), (1,4) \big),
 \quad \big((2,2)\big),\quad \big((3,3)\big) .
\end{alignat*} 
This is unusual in our experience because $\alpha_1=\alpha_2$,
$\alpha_3=\alpha_4$, $\beta_1=\beta_3$ and $\beta_2=\beta_4$.
Compositions are constant after two compositions:
\[ \alpha_1\alpha_3 = 1,\ \alpha_1^2 = 2,\ \alpha_3^2=3, \AND \alpha_3\alpha_1 = 4 \]
and
\[ \beta_1\beta_2 = 1,\ \beta_2^2 = 2,\ \beta_1^2 = 3, \AND \beta_2\beta_1 = 4 .\]
It is $(2,-2)$-periodic.

\end{eg}

\begin{eg}
This final example gives some variants on Example~\ref{twelve}.
For any $m\ge4$, consider the $m\times m$ example which consists of all flips
$\big((i,j),(j,i)\big)$ and fixed points $\big((i,i)\big)$ except when
exactly one of $i,j$ belongs to $\{1,m\}$.  These belong to the $4(m-2)$-cycle
\begin{align*}
 \big( &(2,1),(1,2),(3,1),(1,3),\dots,(m-1,1),(1,m-1),\\
 &(m,2),(2,m),(m,3),(3,m),\dots,(m,m-1),(m-1,m) \big) .
\end{align*}
The maps $\alpha_i$ and $\beta_j$ become constant after three compositions.  

Computer tests show that when $m=2k+2$ is even, the algebra 
has $(12k,-12k)$-periodicity for $1 \le k \le 9$.
Simple examples show that they are not $(6k,-6k)$ or $(4k,-4k)$ periodic.
Except for $k=1$ (Example~\ref{twelve}) in which an exhaustive
computer check was performed, the computer tested a 
random set of a million words of length $12k$ and found the
algebras to be $(12k,-12k)$-periodic.  Experience shows that 
failure of periodicity exhibits itself within a small number of examples.

On the other hand, when $m$ is odd, these examples are aperiodic.  
Since $e_i$ commutes with $f_i$, one sees that $e_i^k$ commutes with $f_i^k$.
Therefore the bijection $\gamma$ of $m^k$ demonstrating periodicity must
map $e_i^k$ to $f_i^k$.  Hence if the algebra is $(k,-k)$-periodic,
one would need to have the identity $e_1^k f_2^k = f_1^k e_2^k$.

For the $5\times 5$ $12$-cycle algebra, we find by induction that
\begin{alignat*}{2}
e_1^{2k+1} f_2^{2k+1} &= f_1^{k+1} f_2 f_5^{k-1} e_5 e_2^{2k} &\qfor& k\ge1 \\
e_1^{2k} f_2^{2k}     &= f_1^{k+1} f_4 f_5^{k-2} e_5 e_2^{2k-1}   &\qfor& k \ge 2 .
\end{alignat*}
Hence it is aperiodic.
Similarly for $m=2s+1$, one can show that 
\[
  e_1^k f_2^k = f_1^{k-l-1} f_{2j} f_m^l e_m e_2^{k-1}
\]
for $k \ge m-2$, where $k \equiv j-2 \pmod{s}$ and $l = \lfloor k/s \rfloor - 2$.
So again these algebras are all aperiodic.
\end{eg}

\section{Periodicity and the Structure of $\ca(\Fth)$}
\label{S:periodStructure}

We first provide a different proof of the Kumjian--Pask result that
aperiodicity implies simplicity of $\ca(\Fth)$.
We have already observed that there is a faithful expectation $\Phi$
onto a $(mn)^\infty$-UHF algebra $\fF$.
If we can show that any ideal $\J$ of $\ca(\Fth)$ is
mapped by $\Phi$ into $\J \cap \fF$, then the simplicity
of $\fF$ will imply that $\ca(\Fth)$ is also simple.
To do this, we copy an argument that works for the
Cuntz algebra \cite{Cun} (see also \cite[Theorem~V.4.6]{KRDCalg}).
We show that in the aperiodic case, the canonical expectation
onto the UHF subalgebra is approximately inner.

\begin{thm}\label{ApproxInner}
Let $\Fth$ be aperiodic.
There is a sequence of isometries $W_k \in \ca(\Fth)$
so that
\[ \Phi(A) = \lim_{k\to\infty} W_k^* A W_k \qforal A \in \ca(\Fth) .\]
\end{thm}

\begin{proof}
It suffices to prove the claim for elements of the form $uv^*$
where $u,v \in \Fth$.
Recall that $\Phi(uv^*) = uv^*$ if $d(u)=d(v)$,
and $\Phi(uv^*) = 0$ otherwise.
Moreover in the first case, it suffices to suppose
that $d(u)=d(v)=(k,k)$; for if $d(u)=(k_1,k_2) \le (k,k)$,
then
\[ uv^* = \sum_{d(x)=(k-k_1,k-k_2)} (ux) (vx)^* \]
is the sum of words with the desired degree.

Let $\tau$ be an aperiodic tail constructed by
Theorem~\ref{aperiodic}.
Then there is a finite segment, say $\tau_k$, which has no
$(a,b)$-periodicity for $|a| \le m^k$ and $|b| \le n^k$.
Let $\S_k = \{ x\in \Fth : d(x) = (k,k) \}$.
Then set
\[ W_k = \sum_{x \in \S_k} x \tau_k x^* .\]

Suppose that $d(u)=d(v) = (k,k)$.
Then
\begin{align*}
 W_k^* uv^* W_k &=
  \sum_{x \in \S_k}  \sum_{y \in \S_k} x \tau_k^* x^* u v^* y \tau_k y^*\\
  &= u \tau_k^* \tau_k v^* = uv^*.
\end{align*}

On the other hand, suppose that
\[
 d(u) \vee d(v) \le (k,k) \qand d(v) - d(u) = (a,b) \ne (0,0) .
\]
Then $x^* u v^* y$ will either be 0 or will have the form
$x_0^*y_0$ in reduced form of total degree $(a,b)$.
Therefore $\tau_k^* x_0^*y_0 \tau_k = 0$ because
$\tau_k$ does not have $(a,b)$ periodicity.
So an examination of the calculation above yields
$W_k^* uv^* W_k = 0$.
\end{proof}

\begin{cor}\label{simple}
If $\Fth$ is aperiodic, then $\ca(\Fth)$ is simple.
\end{cor}

\begin{proof}
If $\J$ is a non-zero ideal in $\ca(\Fth)$, let $A$ be a
non-zero positive element.
Then $\Phi(A) = \lim_{k\to\infty} W_k^* A W_k$
belongs to $\J \cap \fF$.
Since $\Phi$ is faithful, $\Phi(A) \ne 0$.
It follows that $\J$ contains the ideal of $\fF$
generated by $\Phi(A)$, which contains the identity
because $\fF$ is simple.
Therefore $\J = \ca(\Fth)$.
\end{proof}

\subsection{Periodic algebras}
We now turn to the structure of $\ca(\Fth)$ when $\Fth$
is periodic.
Our goal is to show that
$\ca(\Fth) \simeq \rC(\bT) \otimes \fA$
where $\fA$ is simple.

Assume that the minimum period is $(a,-b)$ for $a,b > 0$.
By Theorem~\ref{T:periodic}, there is a bijection
$\gamma:\bm^a \to \bn^b$ so that
\[
 e_uf_v = f_{\gamma(u)} e_{\gamma^{-1}(v)}
 \qforal u \in \bm^a \AND v \in \bn^b .
\]

\begin{lem}\label{L:unitary}
Let $\Fth$ be periodic with minimal period $(a,-b)$ and
define $W := \sum_{u \in \bm^a} f_{\gamma(u)} e_u^*$.
Then $W$ is a unitary in the centre of $\ca(\Fth)$.
\end{lem}

\begin{proof}
It is clear that $W$ is unitary since $\{e_u : u \in \bm^a\}$
is a set of Cuntz isometries, as are $\{f_v : v \in \bn^b\}$.
By Proposition~\ref{periodic identity},
\begin{align*}
 e_i W &= \sum_{u \in \bm^a} e_i f_{\gamma(u)} e_u^*
    = \sum_{u' \in \bm^{a-1}\!\!\!,\, j \in \bm}  e_i f_{\gamma(u'j)} e_{u'j}^* \\
 &=  \sum_{u' \in \bm^{a-1}} f_{\gamma(iu')} \sum_{j=1}^m e_j e_j^*\ e_{u'}^*
    = \sum_{u' \in \bm^{a-1}} f_{\gamma(iu')} e_{u'}^* .
\end{align*}
Thus we compute
\begin{align*}
 W e_i &= \sum_{u \in \bm^a} f_{\gamma(u)} e_u^*  e_i
   = \sum_{u' \in \bm^{a-1}\!\!\!,\, j \in \bm}  f_{\gamma(ju')} e_{ju'}^* e_i\\
 &= \sum_{u' \in \bm^{a-1}\!\!\!,\, j \in \bm}  f_{\gamma(ju')} e_{u'}^*(e_j^*e_i)
   = \sum_{u' \in \bm^{a-1}} f_{\gamma(iu')} e_{u'}^*
   = e_i W .
\end{align*}
Similarly, $W$ commutes with each $f_j$;
and hence it lies in the centre of $\ca(\Fth)$.
\end{proof}

\begin{cor} \label{C:f=ew}
Let $\Fth$ be periodic with minimal period $(a,-b)$.
Then $f_{\gamma(u)} = e_u W$ for all $u \in \bm^a$.
Also if $u \in \bm^a$ and $v \in \bn^b$, then
$e_u^* f_v = \delta_{u,\gamma^{-1}(v)} W$.
\end{cor}

\begin{proof}
\[
 e_u W = W e_u =
 \sum_{u'\in \bm^a} f_{\gamma(u')} (e_{u'}^* e_u) = f_{\gamma(u)} .
\]
Therefore
\[ e_u^* f_{\gamma(u')} = e_u^* e_{u'} W = \delta_{u,u'} W . \qedhere\]
\end{proof}

\subsection{The direct integral decomposition}
Let $\lambda_\tau$ be the representation
constructed in Corollary~\ref{C:periodic}.
This representation is a faithful representation of $\ca(\Fth)$
by \cite{DPYdiln}.
By construction, the symmetry group $H_\tau = \bZ(a,-b)$.
By \cite{DPYatomic}, $\lambda_\tau$ decomposes as a
direct integral of irreducible representations of type 3bii.
It will be helpful to see how this is done in this case
using the extra structure in our possession.

Following the explicit construction of Example~\ref{3a_repn},
we have described $\lambda_\tau$ as an inductive limit
of copies of the left regular representation.
Indeed, $\Fth$ acts on the set
$\G_\tau =  \dlim_{\rightarrow}  (\G_s, \rho_s)$
where the $\rho_s$ are injections of $\G_s$ into $\G_{s+1}$
determined by the word $\tau$, and $\iota_s$ are injections
of $\G_s$ into $\G_\tau$.
We formed $\H_\tau = \ltwo(\G_\tau)$ and obtained a
faithful representation of $\ca(\Fth)$ via the action
\[
 \lambda_\tau(w) \xi_{\iota_s(x)} = \xi_{\iota_s(wx)}
 \qforal w,x \in \Fth \AND s \ge 0.
\]

To understand the way that $W$ of Lemma~\ref{L:unitary}
acts on a basis vector $\xi$, observe that there is a unique
word $u \in \bm^a$ with $\xi \in \ran \lambda_\tau(e_u)$.
The unitary $W$ acts on $\xi$ by pulling back $a$ steps
along the blue edges to $\lambda_\tau(e_u)^* \xi$,
and then pushing forward via $\lambda_\tau(f_{\gamma(u)})$.
This can be computed at a basis vector $\xi$
by representing it as $\xi_{\iota_s(x)}$ with $d(x) \ge (a,0)$
by choosing $s$ sufficiently large.
Then $x$ factors as $x = e_u x'$ with $|u|=a$, and
$\lambda_\tau(W) \xi_{\iota_s(x)} =  \xi_{\iota_s(f_{\gamma(u)}x')}$.

Put an equivalence relation $\sim$ on $\G_\tau$ by taking
the equivalence classes to be the orbits of each basis vector
under powers of $W$.  That is, $x \sim y$ if and only if there
is an integer $k \in \bZ$ so that $W^k \xi_x = \xi_y$.
Let the equivalence classes be denoted by $[x]$, and
set $\W_{[x]} = \spn\{ \xi_y : y \in [x] \}$.
Identify each $\W_{[x]}$ with $\ltwo(\bZ)$ by fixing
a representative $x' \in [x]$ and sending the
standard basis $\{ \eta_k : k \in \bZ\}$ for $\ltwo(\bZ)$
to $\W_{[x]}$ by setting $J_{[x]}\eta_k = W^k \xi_{x'}$.

We wish to choose the representative for each
equivalence class in a consistent way.
So begin with the element $x_0 = \iota_0(\mt)$.
Let
\[
 \S = \{ x \in \G_\tau :
 \xi_x = \lambda_\tau(e_u e_v^* f_w^*) \xi_{x_0}
 \FOR u,v \in \Fm,\ w \in \Fn,\ |w| < b . \}
\]
We claim that $\S$ intersects each equivalence
class of $\G_\tau/\!\!\sim$ in a single point.

Let $x \in \G_\tau$, and choose $s$ so that $\iota_s(\G_s)$
contains $x$.
Then with $\xi = \xi_{\iota_s(\mt)}$, we have words
$e_uf_v$ and $f_{w}e_{u'}$ so that
$\lambda_\tau(e_uf_v) \xi = \xi_x$ and
$\lambda_\tau(f_w e_{u'}) \xi = \xi_{x_0}$.
If $|v| \equiv -k \pmod{b}$ for $1 \le k < b$,
we replace $\xi$ by $\xi_{\iota_{s+k}(\mt)}$ so
that $|v|$ is a multiple of $b$.
Factor $f_v = f_{v_1} \dots f_{v_s}$ where $|v_i|=b$
and factor $f_w = f_{w_0} f_{w_1} \dots f_{w_t}$ where
$0 \le |w_0| < b$ and $|w_j|=b$ for $1 \le j \le t$.
Then set $u_i = \gamma^{-1}(v_i)$ and $x_i = \gamma^{-1}(w_i)$
and use Corollary~\ref{C:f=ew}:
\begin{align*}
 \xi_x &= \lambda_\tau(e_uf_ve_{u'}^*f_w^*) \xi_{x_0} \\&=
 \lambda_\tau(e_u) \lambda_\tau(f_{v_1} \dots f_{v_s})
 \lambda_\tau(e_{u'}^*) \lambda_\tau(f_{w_t}^* \dots f_{w_1}^*)
 \lambda_\tau(f_{w_0})^* \xi_{x_0} \\&= \lambda_\tau(e_u)
 \lambda_\tau(W^s e_{u_1} \dots e_{u_s}) \lambda_\tau(e_{u'}^*)
 \lambda_\tau(W^{t*} e_{x_t}^* \dots e_{x_1}^*)
 \lambda_\tau(f_{w_0}^*) \xi_{x_0} \\&=
 \lambda_\tau(W^{s-t}) \lambda_\tau(e_{u''}e_{x''}^* f_{w_0}^*) \xi_{x_0}
 =:  \lambda_\tau(W^{s-t}) \xi_y .
\end{align*}
where $u'' = u u_1 \dots u_s$ and $x''=x_1 \dots x_t u'$
and $\xi_y = \lambda_\tau(e_{u''}e_{x''}^* f_{w_0}^*) \xi_{x_0}$.
Therefore $[x] \cap \S$ contains $y$.

Uniqueness follows because there is an essentially unique
way to write $\xi_x = \lambda_\tau(e_uf_ve_{u'}^*f_w^*) \xi_{x_0}$
except for reducing the word because of redundancies.
As
\[
 \xi_x = \lambda_\tau(W^k) \xi_y =
 \lambda_\tau(W^k e_{u''}e_{v''}^*f_{w_0}^*) \xi_{x_0} ,
\]
one sees $W^k e_{u''}e_{v''}^*f_{w_0}^*$ has
degree $(|u''|-|v''|-ka, kb-|w_0|)$.  Since $kb-|w_0|$
is not in $[1-b,0]$, this point does not lie in $\S$.

The fact that $W$ lies in the centre of $\ca(\Fth)$
means that for each $w\in\Fth$ and $\lambda_\tau(w) \xi_x = \xi_y$,
then $\lambda_\tau(w W^k) \xi_x = \lambda_\tau(W^k) \xi_y$.
Thus $\lambda_\tau(w)$ maps each subspace $\W_{[x]}$
onto another subspace $\W_{[y]}$.

Let $U$ be the bilateral shift on $\ltwo(\bZ)$.
Then one sees that $\lambda_\tau(W)$ acts as a shift
$J_{[x]} U J_{[x]}^*$ on every $\W_{[x]}$.
Hence $\lambda_\tau(W) \simeq U \otimes I$ is a
bilateral shift of infinite multiplicity.
In particular, the spectrum of $W$ is the whole circle $\bT$;
and the spectral measure of $\lambda_\tau(W)$ is
Lebesgue measure.

Let us consider how $\lambda_\tau(e_i)$ acts on $\W_{[x]}$.
Say $[x]\cap\S = \{y\}$ where
$\xi_y = \lambda_\tau(e_{u}e_{v}^* f_{w}^*) \xi_{x_0}$.
Then $[e_ix]\cap\S$ is $z$ where
\[ \xi_z = \lambda_\tau(e_{iu}e_{v}^* f_{w}^*) \xi_{x_0} .\]
This is in reduced form unless $u=\mt$ and $v= v'i$, in which
case, after cancellation,
$\xi_z = \lambda_\tau(e_{v'}^* f_{w}^*) \xi_{x_0}$.
Either way, we see that $z \in \S$.
Hence
\[
 \lambda_\tau(e_i)|_{\W_{[x]}} = J_{[z]}J_{[y]}^* = J_{[e_ix]}J_{[x]}^* .
\]

Similarly we analyze $\lambda_\tau(f_j)$.
This comes down to understanding the representative
$[f_jx] \cap \S$.
Again we use the representative $y \in [x]$ with
$\xi_y = \lambda_\tau(e_{u}e_{v}^* f_{w}^*) \xi_{x_0}$.
If $|w| \ge 1$, write $w = w'i$.  Then there is a unique
word $\tilde v \in \Fn$ with $|\tilde v|=b-1$ so that
$ \lambda_\tau(e_{v}^* f_{w}^*) \xi_{x_0} =
\lambda_\tau(f_{\tilde v}f_{\tilde v}^* e_v^*f_{w}^*) \xi_{x_0}$
is in the range of $\lambda_\tau(f_{\tilde v})$.
Therefore using the commutation relations
$f_je_u = e_{u'}f_{j'}$ and $e_v^*f_i^* = f_{i'}^* e_{v'}^*$,
\begin{align*}
 \lambda_\tau(f_j e_u e_v^* f_w^*) \xi_0 &=
 \lambda_\tau(f_j e_u f_{\tilde v}f_{\tilde v}^* e_v^* f_w^*) \xi_0 \\ &=
 \lambda_\tau(e_{u'} f_{j' \tilde v} f_{i'\tilde v}^* e_{v'}^* f_{w'} ) \xi_0 \\ &=
 \lambda_\tau(e_{u'} e_{j' \tilde v}WW^* e_{i'\tilde v}^* e_{v'}^* f_{w'} ) \xi_0
 \\&= \lambda_\tau(e_{u'j' \tilde v}e_{v'i'\tilde v}^* f_{w'} ) \xi_0
\end{align*}
So again, we see that the canonical representative $y$ in $[x]$ is
carried by $f_j$ to $f_jy$, the canonical representative in $[f_jx]$;
whence $\lambda_\tau(f_j)|_{\W_{[x]}} = J_{[f_jx]} J_{[x]}^*$.

However, when $w=\mt$, one obtains
\begin{align*}
 \lambda_\tau(f_j e_u e_v^*) \xi_0 &=
 \lambda_\tau(e_{u'} f_{j' \tilde v} f_{\tilde v}^* e_{v'}^*) \xi_0 \\ &=
 \lambda_\tau(e_{u'} e_{j' \tilde v}W f_{\tilde v}^* e_{v'}^*) \xi_0
 \\&= \lambda_\tau(W e_{u'j' \tilde v}e_{v''}^* f_{w'} ) \xi_0
\end{align*}
where $|w'| = |\tilde v| = b-1$.
So when $\lambda_\tau(f_j)$ maps $\W_{[x]}$ to $\W_{[f_jx]}$,
it is acting like the bilateral shift, namely
$\lambda_\tau(f_j)|_{\W_{[x]}} = J_{[f_jx]} U J_{[x]}^*$.

Form a Hilbert space $\K = \ltwo(\G_\tau/\!\!\sim)$
with basis $\{ \zeta_{[x]} ; [x] \in \G_\tau/\!\!\sim \}$.
One can see that $\lambda_\tau$ is unitarily equivalent
to a representation $\pi_\tau$ on
$\ltwo(\G_\tau/\!\!\sim) \otimes L^2(\bT)$ given by
\begin{align*}
 \pi_\tau(e_i) \zeta_{[x]} \otimes h &= \zeta_{[e_i x]} \otimes h \\
 \pi_\tau(f_j) \zeta_{[x]} \otimes h &=
 \begin{cases}
  \zeta_{[f_j x]} \otimes h &\qif [x]=[e_u e_v^* f_w^* x_0],\ 1 \le |w| < b \\
  \zeta_{[f_j x]} \otimes zh &\qif [x]=[e_u e_v^* x_0]
 \end{cases}
\end{align*}
for all $[x]\in\G_\tau/\!\!\sim$ and $h \in L^2(\bT)$.

Define a representation $\sigma_1$ by
\[
 \sigma_1(w) \zeta_{[x]} = \zeta_{[y]}
 \qif \lambda_\tau(w) \W_{[x]} = \W_{[y]} .
\]
i.e.\ $\sigma_1(w) \zeta_{[x]} = \zeta_{[wx]}$.
Then for each $z\in\bT$, define
\[
 \sigma_z(e_i) = \sigma_1(e_i) \AND
 \sigma_z(f_j) = z \sigma_1(f_j)
\]
and extend to $\Fth$.
It is not difficult to see that $\sigma_z$ is unitarily
equivalent to another atomic representation $\rho_{z^b}$ given by
\begin{align*}
 \rho_{z^b}(e_i) \zeta_{[x]} &= \zeta_{[e_i x]}  \\
 \rho_{z^b}(f_j) \zeta_{[x]} &=
 \begin{cases}
  \zeta_{[f_j x]} &\qif [x]=[e_u e_v^* f_w^* x_0],\ 1 \le |w| < b \\
  z^b \zeta_{[f_j x]} &\qif [x]=[e_u e_v^* x_0]
 \end{cases}
\end{align*}
In particular, $\sigma_z \simeq \sigma_w$ if
and only if $z^b=w^b$.

The representations $\sigma_z$ are irreducible by
\cite[Lemma~8.14]{DPYatomic}
because the symmetry group $H_{\sigma_z} = \{0\}$
is the trivial subgroup of $\bZ^2/\bZ(a,-b)$.
Note that $\sigma_z(W) = z^b I$.

  From this picture, one can see how to decompose the
representation $\lambda_\tau$ as a direct integral.
Indeed,
\[
 \lambda_\tau \simeq \int_0^{2\pi/b\ \oplus} \sigma_{e^{it}} \,dt
 \simeq \int_0^{2\pi\ \oplus} \rho_{e^{it}} \,dt.
\]

In particular, notice that the C*-algebra $\fA = \sigma_z(\ca(\Fth))$
is independent of $z$.

\subsection{An expectation}
We need to build a somewhat different expectation
in the periodic case.

\begin{thm}\label{T:Psi}
If $\Fth$ is a periodic semigroup with minimal period $(a,-b)$,
define
\[ \Psi(X) = \int_\bT \gamma_{z^b,z^a}(X) \,dz . \]
Then $\Psi$ is a faithful, approximately inner expectation
onto
\[ \ca(\fF,W) \simeq \rC(\bT) \otimes \fF \simeq \rC(\bT,\fF) .\]
\end{thm}

\begin{proof}
As an integral of automorphisms, $\Psi$ is evidently a
faithful completely positive map.
Suppose that $w = e_uf_ve_{u'}^*f_{v'}^*$ is a word of degree
$(k,l)$ where $k = |u|-|u'|$ and $l=|v|-|v'|$.
Then
\[
 \Psi(e_uf_ve_{u'}^*f_{v'}^*) = \int_\bT z^{kb+la}e_uf_ve_{u'}^*f_{v'}^* \,dz
 = \begin{cases} e_uf_ve_{u'}^*f_{v'}^* &\qif kb+la=0\\
 0 &\quad\text{otherwise.}\end{cases}
\]
That is, $\Psi(w) = w$ if $d(w) \in \bZ(a,-b)$ and is $0$ otherwise.
Therefore this is an idempotent map, and so it is an expectation.
Since the degree map is a homomorphism, the range of $\Psi$
is a C*-subalgebra of $\ca(\Fth)$.

The range contains $\fF$ as this is spanned by words
of degree $(0,0)$ and $W$, which has degree $(-a,b)$.
The typical word of degree $(-pa,pb)$ is
$w = e_uf_ve_{u'}^*f_{v'}^*$ where
$|u|-|u'|= -pa$ and $|v|-|v'|=pb$.
For convenience, consider $p \ge0$.
If $|u| \equiv -k \pmod{a}$ for $0 \le k < a$ and
$|v| \equiv -l \pmod{b}$ for $0 \le l < b$, then
\[ w = \sum_{x \in \Fth, d(x)=(k,l)} e_uf_v x x^*e_{u'}^*f_{v'}^* .\]
Hence $w$ is in the span of words of the form
$e_uf_ve_{u'}^*f_{v'}^*$ for which
\[
 |u| \equiv |u'| \equiv 0 \pmod a \qand
 |v| \equiv |v'| \equiv 0 \pmod b .
\]
For such words, we may split each word $u,u'$ into words
of length $a$ and each $v,v'$ into words of length b as
\begin{align*}
  w' &= e_{u_1\dots u_s}f_{v_1\dots v_{t+p}}
  e_{u'_1\dots u'_{s+p}}^* f_{v'_1 \dots v'_t}^* \\&=
  W^p e_{u_1\dots u_s \gamma^{-1}(v_1) \dots \gamma^{-1}(v_{t+p})}
  e_{\gamma^{-1}(v'_1) \dots \gamma^{-1}(v'_t) u'_1\dots u'_{s+p} }^*
\end{align*}
Therefore all of these words lie in $\ca(\Fth,W)$.
So this identifies the range of $\Psi$ as $\ca(\fF,W)$.

Now $W$ lies in the centre of $\ca(\Fth)$; and it is easy to
check that $\ca(W) \cap \fF = \bC I$.
So a dense subalgebra of $\ca(\fF,W)$ is given by
the polynomials $\sum_k F_kW^k$ where the sum is finite
and $F_k \in \fF$.  Since the spectrum of $W$ is $\bT$,
\[
 \big\| \sum_k F_kW^k \big\| =
 \sup_{z \in \bT} \big\| \sum_k F_k z^k \big\| .
\]
Now a routine modification of the proof of F\'ejer's Theorem shows that
$\ca(\fF,W) \simeq \rC(\bT,\fF) \simeq \rC(\bT) \otimes \fF$.

The last part of the proof is to establish that $\Psi$ is
approximately inner.
The argument is a modification of the proof of Theorem~\ref{ApproxInner}.

Let $\tau$ be the infinite tail used above whose only
symmetries are $\bZ(a,-b)$.
Then there is a finite segment, say $\tau_k$, such that
$\tau_k^* u^*v \tau_k = 0$ whenever $u,v \in \Fth$
with
\[
 d(u)-d(v) \not\in \bZ(a,-b) \qand \max\{d(u),d(v)\} \le (ka,kb) .
\]
Let $\S_k = \{ x\in \Fth : d(x) = (ak,bk) \}$.
Then set
\[ W_k = \sum_{x \in \S_k} x \tau_k x^* .\]

Suppose that $w=uv^*$ for $u,v\in\Fth$ such that
$d(w) \in \bZ(a,-b)$ and $\max\{d(u),d(v)\} \le (ka,kb)$.
Then as before, we may write $w$ as a sum of words with the
same property and the additional stipulation that
$d(u) \vee d(v) = (ka,kb)$.
Say, as a typical case, that $d(u)=(ka,(k-p)b)$ and $d(v) = ((k-p)a,kb)$.
Then $uv^* = W^{-p} (u W^p v^*)$, and the term $u W^p v^*$ splits
as a sum of words $xy^*$ with $d(x)=d(y)=(pa,pb)$.
Thus if $d(x_0) = d(y_0) = (pa,pb)$,
\begin{align*}
 W_k^* W^{*p} x_0 y_0^* W_k &=
  W^{*p} \sum_{x \in \S_k}  \sum_{y \in \S_k} x \tau_k^* x^* x_0 y_0^* y \tau_k y^*\\
  &= W^{*p} x_0 \tau_k^* \tau_k y_0^* = W^{*p} x_0 y_0^*.
\end{align*}

On the other hand, suppose that $w=uv^*$ for $u,v \in \Fth$
with $d(w) \not\in \bZ(a,-b)$ and $\max\{d(u),d(v)\} \le (ka,kb)$.
Then $x^* u v^* y$ will either be 0 or will have the form
$x_0^*y_0$ in reduced form, with
\begin{align*}
 d(y_0)-d(x_0) &=
 (d(v)-d(y)) - (d(u)-d(x)) \\&=
 d(v) - d(u) =
 -d(w) \not\in\bZ(a,-b) .
\end{align*}
Therefore $\tau_k^* x_0^*y_0 \tau_k = 0$;
whence $W_k^* uv^* W_k = 0$.
\end{proof}

Now we `evaluate at 1' and obtain a faithful, approximately
inner expectation of $\fA = \sigma_1(\ca(\Fth))$ onto $\fF$.

\begin{thm}\label{T:Psi1}
There is a faithful, approximately inner expectation $\Psi_1$ of
$\fA$ onto $\fF$ so that the following diagram commutes,
where $\ep_1$ is evaluation at $1$:
\[
 \xymatrix{ \ca(\Fth) \ar[r]^(.6){\sigma_1} \ar[d]_{\Psi} & \fA \ar[d]^{\Psi_1}\\
 \rC(\bT,\fF) \ar[r]^(.6){\ep_1}&\fF }
\]
\end{thm}

\begin{proof}
The first step is to observe that \cite[Lemma~3.4]{DPYdiln}
shows that there is a unitary $U_z := U_{z^b,z^a}$ given by
\[
 U_z \xi_{\iota_s(e_uf_v)} =
 z^{b(|u|-s)+a(|v|-s)} \xi_{\iota_s(e_uf_v)} .
\]
which implements $\gamma_{z^b,z^a}$ in that
\[
 \gamma_{z^b,z^a}(X) = U_z \lambda_\tau(X) U_z^*
 \qfor X \in \ca(\Fth) .
\]
Observe that $U_z$ in constant on the subspaces $\W_{[x]}$,
and thus it determines a unitary $V_z$ on
$\ltwo(\G_\tau/\!\!\sim)$ by $V_z \zeta_{[x]} = t  \zeta_{[x]}$
if $U_z|_{\W_{[x]}} = t P_{\W_{[x]}}$.

Clearly
\begin{align*}
 V_z \sigma_1(X) V_z^* &= \sigma_1( U_z \lambda_\tau(X) U_z^*)
 = \sigma_1( \gamma_{z^b,z^a}(X) ) .
\end{align*}
This determines an automorphism $\psi_z$ of $\fA$.
Define a map
\[ \Psi_1(A) = \int_\bT \psi_z(A) \,dz \qfor A \in \fA .\]
It follows that
\[ \Psi_1(\sigma_1(X)) = \sigma_1(\Psi(X)) \qfor X \in \ca(\Fth) .\]
In particular, $\Psi_1$ is a faithful expectation of $\fA$
onto $\sigma_1(\ca(\fF,W))$.

Now $\sigma_1(W) = I$.
So $\sigma_1(\ca(\fF,W)) = \sigma_1(\fF) \simeq \fF$ because
$\fF$ is simple.
The ideal of $\ca(\fF,W)$ generated by $W-I$ is contained in
$\ker \sigma_1$.
But this ideal is evidently $\rC_0(\bT\setminus\{1\}) \otimes \fF$.
Thus $\sigma_1|_{\ca(\fF,W)} = \ep_1$ is evaluation at $1$.

Finally, since $\Psi$ is approximately inner, we obtain
for any $A = \sigma_1(X)$ in $\fA$ that
\begin{align*}
 \Psi_1(A) &= \sigma_1(\Psi(X)) =
 \lim_{k\to\infty} \sigma_1(W_k^* X W_k) \\&=
 \lim_{k\to\infty} \sigma_1(W_k^*) A \sigma_1(W_k) .
\end{align*}
So $\Psi_1$ is approximately inner.
\end{proof}

\begin{cor}\label{periodic simple}
$\fA$ is simple.
\end{cor}

\begin{proof}
If $\J$ is a non-zero ideal, then it contains a positive element $A$.
Since $\Psi_1$ is faithful, $\Psi_1(A) \ne 0$.
As $\Psi_1$ is approximately inner, $\Psi_1(A)$ belongs to
$\fF \cap \J$.
But $\fF$ is simple and unital, so $\J$ contains $\fF$
and thus is all of $\fA$.
\end{proof}

\begin{cor}\label{periodic algebra}
Let $\Fth$ have minimal period $(a,-b)$.  Then
\[ \ca(\Fth) \simeq \rC(\bT) \otimes \fA .\]
\end{cor}

\begin{proof}
Here it is more convenient to use the representations $\rho_z$.
The analysis above applied to each representation $\rho_z$
shows that $\ker \rho_z = \ip{W - z I}$.
Define a map $\phi$ from $\ca(\Fth)$ to $\rC(\bT,\fA)$
by $\phi(X)(z) = \rho_z(X)$.
Checking continuity is straightforward,
as is surjectivity.
That the map is injective follows from the direct integral
decomposition of the faithful representation $\lambda_\tau$.
\end{proof}

As a corollary, we obtain a special case of the result of
Robertson and Sims \cite{RobSims} on the simplicity of
higher rank graph algebras.

\begin{cor}\label{simplicity}
Let $\Fth$ be a rank 2 graph with a single vertex.
Then $\ca(\Fth)$ is simple if and only if
$\Fth$ is aperiodic.
\end{cor}

\begin{cor}\label{centre}
Let $\Fth$ be a periodic rank 2 graph with a single vertex.
Then the centre of $\ca(\Fth)$ is $\ca(W) \simeq\rC(\bT)$.
\end{cor}



\end{document}